\documentclass[11pt]{amsart}

\usepackage{amscd,latexsym,amsthm,amsfonts,amssymb,amsmath,amsxtra,eucal}
\usepackage[all]{xy}

\usepackage[pdftex,bookmarks,colorlinks,pdfmenubar]{hyperref}

\usepackage{verbatim}
\usepackage{mathabx}
\usepackage{mathrsfs}
\usepackage{bbm}

\renewcommand\theequation{\thesection.\arabic{equation}}

\newcommand{\Qlb}{\mathbb{\bar Q}_\ell}

\newcommand{\Ad}{{\mathrm{Ad}}}

\newcommand{\GL}{{\mathrm{GL}}}

\newcommand{\Hom}{{\mathrm{Hom}}}

\newcommand{\Res}{{\mathrm{Res}}}

\newcommand{\SO}{{\mathrm{SO}}}

\newcommand{\Sp}{{\mathrm{Sp}}}

\newcommand{\Span}{{\mathrm{Span}}}

\newcommand{\wt}{\widetilde}

\newtheorem{thm}{Theorem}[section]
\newtheorem{cor}[thm]{Corollary}

\newtheorem{prop}[thm]{Proposition}
\newtheorem {conj}[thm]{Conjecture}
\newtheorem {ques/conj}[thm]{Question/Conjecture}

\newtheorem{rmk}[thm]{Remark}

\usepackage[numbers]{natbib}

\newcommand{\select}[1]{{\it{#1}}}

\begin{document}
\renewcommand{\theequation}{\arabic{equation}}
\numberwithin{equation}{section}

\title[upper bound of the multiplicities of Fourier-Jacobi models]{On the upper bound of the multiplicities of Fourier-Jacobi models over finite fields}

\date{\today}







\author[Fang Shi]{Fang Shi}

\address{School of Mathematical Sciences, Zhejiang University, Hangzhou 310027, Zhejiang, P.R. China}

\email{11935007@zju.edu.cn}


\maketitle

\begin{abstract}
In general the multiplicity one theorem fails for Fourier-Jacobi models over finite fields. In this paper we prove that there is an upper bound for the multiplicities of Fourier-Jacobi models which is independent of $q$. As a  consequence, we prove  a conjecture of Hiss and Schr\"oer.
\end{abstract}

\section{Introduction and Main Results}
~

\subsection{Conventions}
Let $q$ be a power of a prime $p\neq 2$.
For schemes over $\mathbb F_q$, we denote the Frobenius by $F$.  We fix an algebraic closure $\mathrm k_q$ of $\mathbb F_q$. 
We will abuse the notation by denoting the pullback of $F$ to $\mathrm k_q$ again by $F$. We say that an algebraic group $H_0$ over $\mathbb F_q$ is reductive if its pullback $H$ to $\mathrm k_q$ is reductive (\select{i.e.}, $H$ is smooth connected and   the unipotent radical of $H$ is trivial).
For a reductive group $G$ over $\mathbb F_q$, we denote
the $F$-rank of $G$ by $\sigma(G)$.  We choose a prime number $\ell$ different from $p$.
 The set of (the isomorphism classes of) irreducible $\Qlb$-representations of the finite group $G^F$ is denoted by $\mathcal E(G)$. We implicitly use an identification $\Qlb \simeq \mathbb C$. The results of this paper do not rely on the choice of the identification. When $q$ and $\ell$ are clear, we choose a nontrivial additive character $\psi: \mathbb F_q \to \Qlb^\times$. For a finite group $K$, we will frequently denote $\dim \Hom_K(\sigma,\pi)$ by $\langle \sigma,\pi\rangle_K$ for a pair of $\Qlb$-representations $(\sigma,\pi)$ of $K$. We will abuse the notation by denoting the linear extension of $\langle -,-\rangle_K$ to virtual characters of $K$ again by $\langle-,-\rangle_K$. Let $\mathbb Z_+$ be the set of positive integers.



We often identify an algebraic group $H_0$ over $\mathbb F_q$ with its pullback $H$ to $\mathrm k_q$ equipped with the Frobenius endomorphism $F:H \to H$. This should cause no confusion.
\subsection{Main Results}
In general the multiplicity one theorem fails for Fourier-Jacobi models over finite fields (see Section 4.4 of \cite{LMS}). In this paper, we prove that there is an upper bound for the multiplicities of Fourier-Jacobi models which is independent of $q$, see Theorem \ref{mt}, Theorem \ref{mt2} and Remark \ref{rmmt}. In particular, our results imply Conjecture 3.1 of \cite{HS}.
\begin{conj}[Hiss and  Schr\"oer] \label{conj}
Let $G_n= \mathrm U_n$ or $\Sp_{2n}$.
There is a  function $f:\mathbb Z_+ \to \mathbb Z_+$ that is independent of $q$ satisfying the following
$$ \langle \pi ,\pi' \otimes \omega_{G_n(\mathbb F_q),\psi}\rangle_{G_n(\mathbb F_q)} \leq f (n),$$
for all irreducible representations $\pi,\pi'$ of $G_n(\mathbb F_q)$, where $\omega_{G_n(\mathbb F_q),\psi}$ is the Weil representation of $G_n(\mathbb F_q)$ corresponding to $\psi$.
\end{conj}
From now on, we denote the Weil representation of $\Sp_{2n}(\mathbb F_q)$ corresponding to $\psi: \mathbb F_q \to \Qlb^\times$ by $\omega_\psi$. The dual representation of $\omega_\psi$ is denoted by $\omega_\psi^\vee$.
In this paper we will prove the following. 
\begin{thm}\label{pmt} There is a function $f: \mathbb Z_+ \to \mathbb Z_+$ which is independent of q, so that the following holds: For any $\pi,\pi'\in \mathcal E(\Sp_{2n})$, we have
$$\dim \Hom_{\Sp_{2n}(\mathbb F_q)} (\pi \otimes \omega_\psi^\vee ,\pi') \leq f(n).
$$
\end{thm}
The above theorem can be viewed as an estimation for the basic case of Fourier-Jacobi models. We will extend the above to a result concerning general Fourier-Jacobi models in Theorem \ref{mt2}.
\begin{rmk}
The case for unitary groups in Conjecture \ref{conj} is proved in Theorem 4.5 of \cite{LMS}.
\end{rmk}
~

\subsection{An overview}
In this subsection we briefly introduce the strategy. We fix a connected reductive group $G$ over $\mathbb Z$ in this subsection. We denote the set (of the isomorphism  classes) of irreducible $\mathbb C$-representations of $G(\mathbb F_q)$ by $\mathcal E(G(\mathbb F_q))$.

 Let $r=1$ or $r=2$.
 Let $\{d_q^\circ(.)\}_q$ be a family of functions 
$$
d_q^\circ:  (\mathcal E(G(\mathbb F_q)))^r \to \mathbb N,   
$$
where  we set $(\mathcal E(G(\mathbb F_q)))^1$ to be $\mathcal E(G(\mathbb F_q))$ and set $(\mathcal E(G(\mathbb F_q)))^2$ to be $\mathcal E(G(\mathbb F_q))\times \mathcal E(G(\mathbb F_q))$.
We extend each $d_q^\circ$ linearly to 
 virtual characters and get the resulting family $\{d_q (.)\}_q$. 
A typical problem is that whether there is an upper bound for $\{d_q^\circ(.)\}_q$, \select{i.e.}, whether the following is finite:
$$
\sup_q \max_{\chi \in (\mathcal E(G(\mathbb F_q)))^r} d_q^\circ (\chi).
$$
We find it possible to answer this question in two steps.

First of all, we deal with the problem of whether there is an upper bound for 
$d_q(.)$ restricted to 
Deligne-Lusztig characters. Such a problem can sometimes be easier to tackle using some geometry methods. If the answer is affirmative, we denote an upper bound by $A$. In this paper  the above step corresponds to Theorem \ref{MF} and Corollary \ref{bu}.

 We proceed by realizing arbitrary $\iota \in \mathcal E(G(\mathbb F_q))$ as a subrepresentation of  a representation $\rho_\iota$ of $G(\mathbb F_q)$ whose character is a linear combination of Deligne-Lusztig characters. Suppose we can write the character of $\rho_\iota$ as $\sum\limits_i c_i^\iota R_i$ where $R_i$ are Deligne-Lusztig characters and $\sum\limits_i |c_i^\iota|$ are bounded for $\iota$ varying in irreducible representations of $G(\mathbb F_q)$ for all $q$. We denote an upper bound of $\sum\limits_i |c_i^\iota|$ by $B$. This step corresponds to Proposition \ref{embed}.

 Combining the above two, we easily find that $A B^r$ is a desired upper bound for the values of the functions $d^\circ_q$ which is independent of $q$.

A similar argument is adopted in \cite{S}, where $r=1$ and $d^\circ_q$  is the dimension of the maximal $H(\mathbb F_q)$-invariant space for a fixed spherical subgroup $H$ of $G$.

\subsection*{Acknowledgement}
The author would like to thank Dongwen Liu and 	Zhicheng Wang for their helpful discussion and suggestion. The author is also grateful to the referee for careful reading and useful comments.

\section{Preliminaries}
In this section we collect the ingredients that we need.
~

\subsection{Around the Deligne-Lusztig characters}\label{DLc}
Fix a prime $\ell$ different from $p$. For a connected reductive group $G$ over $\mathbb F_q$, the virtual character $R^G_{T,\theta}$ of $G^F$ is defined in \cite{DL}. 
 Here $T$ is an $F$-stable maximal torus in $G$ and $\theta: T^F \to  \Qlb^\times$ is a character.
We say that a class function $f:G^F \to  \Qlb$ is uniform if $f$ is in the $\Qlb$-span of Deligne-Lusztig characters.  
\begin{rmk} \label{rell}
The Deligne-Lusztig characters $R^G_{T,1}$ is independent of the choice of $\ell$, since it is integer-valued. See Proposition 3.3 and the beginning of Section 4 of \cite{DL}.
The virtual character $R^G_{T,\theta}$ only depends on $\theta: T^F \to \Qlb^\times$.   
Indeed, if we fix an inclusion $i_\ell:\bar {\mathbb Q} \hookrightarrow  \Qlb$ for all primes $\ell \neq p$ and view $\theta$ as a composition of  $\theta_0: T^F \to  \bar{\mathbb Q}^\times $ with $i_\ell$, then we see different choices of $\ell$ give the same theory of Deligne-Lusztig characters by the character formula Theorem 4.2 of  \cite{DL}.
\end{rmk}
In this subsection we recall some definitions concerning the theory of Deligne-Lusztig characters; see \cite{L1} or Chapter 11.5 of \cite{DM}.


 Let $G^*$ be the dual group of $G$. There is a natural bijection between the set of $G^F$-conjugacy classes of $(T,\theta)$ and the set of $G^{*F}$-conjugacy classes of $(T^*,s)$, where $T^*$ is an $F$-stable maximal torus in $G^*$ and $s\in T^{*F}$.  We will also denote $R^G_{T,\theta}$ by $R^G_{T^*,s}$ if $(T,\theta)$ corresponds to $(T^*,s)$.

 For a semisimple element $s\in G^{*F}$, define its Lusztig series as 
$$
\mathcal E(G,s)=\{\pi \in \mathcal E(G): \langle \pi ,R^G_{T^*,s} \rangle \neq 0 ~for~some~T^*~containing~s\}.
$$
Then we have 
$$\mathcal E(G) = \bigsqcup_{(s)} \mathcal E(G,s)
$$
where $(s)$ runs over the (rational) conjugacy classes of semisimple elements.
There is a bijection
$$\mathfrak L_s: \mathcal E(G,s) \to \mathcal E(C_{G^*}(s),1),
$$
extended by linearity to a map between virtual characters satisfying that
$$\mathfrak L_s((-1)^{\sigma(G)} R^G_{T^*,s})=(-1)^{\sigma (C_{G^*}(s))}R^{C_{G^*}(s)}_{T^*,1}.
$$
\begin{rmk}
The definition of Deligne-Lusztig characters for a disconnected group $G$ is used in the above bijection. If $G$ is a disconnected group whose identity component $G^\circ$ is reductive, we write $R_{T,\theta}^G:= Ind^{G^F}_{G^{\circ,F}} R_{T,\theta}^{G^\circ}$, where $T$ is a rational maximal torus in $G^\circ$ and $\theta: T^F \to \Qlb^\times$ is a character. 
 See \select{e.g.}, Chapter 11.5 of \cite{DM} for details.
\end{rmk}
~

\subsection{Irreducible representations of $\GL_n$ and $\mathrm{U}_n$}
In this subsection we assume that $G=\GL_n$ or $\mathrm U_n$. We see from Section 3 of \cite{LS} the following.

\begin{prop} \label{Airep}
Every irreducible representation of $G^F$ is of the form
$$
R(G,L,\theta,V)=(-1)^{\sigma(G)-\sigma(L)}|W_L|^{-1} \sum\limits_{w\in W_L} Tr(w\wt F,V) R^G_{T_w,\theta_w}
$$
where
\begin{itemize}
\item $L$ is an $F$-stable reductive connected closed subgroup of $G$, which arises as the connected centralizer of a torus in $G$;
\item $\theta$ is a homomorphism $L^F \to \Qlb^\times$ with certain properties;
\item $\theta_w$ is the restriction of $\theta $ to $T_w$, where $T_w$ is the $F$-stable maximal  torus of $L$ corresponding to $w$ in the sense of Corollary 1.14 of \cite{DL};
\item $W_L$ is the Weyl group of $L$;
\item $V$ is an irreducible representation of $W_L$, and the action is given by the homomorphism $\phi:W_L \to \GL(V)$;
\item There is $w_0\in W_L$, so that $\wt F= s\cdot \phi(w_0) $ as automorphisms of $V$, where $s$ is a root of $1$;
\item For $w\in W_L$, the endomorphism $w \wt F$ is the composition of the endomorphisms $\wt F$ and  $\phi(w)$ of $V$.
\end{itemize}
\end{prop}

For future use, we record the following corollary:

\begin{cor} \label{Airep2}
Let $H$ be a product of general linear groups and unitary groups. Then every irreducible representation of $H^F$ is of the form
$$
R(H,M,\chi,V)=(-1)^{\sigma(H)-\sigma(M)}|W_M|^{-1} \sum\limits_{w\in W_M} Tr(w\wt F,V) R^H_{T_w,\chi_w}
$$
where
\begin{itemize}
\item $M$ is an $F$-stable reductive connected closed subgroup of $H$, which arises as the connected centralizer of a torus in $G$;
\item $\chi$ is a homomorphism $M^F \to \Qlb^\times$ with certain properties;
\item $\chi_w$ is the restriction of $\chi $ to $T_w$, where $T_w$ is a $F$-stable maximal torus of $L$ corresponding to $w$ in the sense of  Corollary 1.14 of \cite{DL};
\item $W_M$ is the Weyl group of $M$;
\item $V$ is an irreducible representation of $W_M$, and the action is given by the homomorphism $\phi: W_M \to \GL(V)$;
\item There is $w_0\in W_M$, so that $\wt F= s\cdot \phi(w_0 )$ as automorphisms of $V$, where $s$ is a root of $1$;
\item For $w\in W_M$, the endomorphism $w \wt F$ is the composition of the endomorphisms $\wt F$ and $\phi(w)$ of $V$.
\end{itemize}

\end{cor}
\begin{rmk}\label{Airep3}
When we take $\pi \in \mathcal E(H,1)$ for $H$ being as in Corollary \ref{Airep2},  we have $\chi_w =1$ for $\chi_w$ on the right hand side of the equation in Corollary \ref{Airep2}. This can be seen from Section 2 of \cite{LS}, or projecting both sides of the above equation to the unipotent part. 
\end{rmk}
~

\subsection{Unipotent irreducible representations of $\SO_{2n+1}$ and $\Sp_{2n}$}
There is no analog of Corollary \ref{Airep2} when $G=\SO_{2n+1}$. For our purpose, the following proposition is adequate.
\begin{prop}\label{Birep}Assume that $q\geq 2^{2n}$.
Let $\pi \in \mathcal E(\SO_{2n+1},1)$. Then $\pi$ is a subrepresentation of a representation $\phi$ whose character $\chi_\phi$ is of the form 
$$\chi_\phi=\sum_i c_i R_i,
$$
where each $R_i$ is of the form $R_{T,1}^{\SO_{2n+1}}$ for some $F$-stable maximal torus $T$ in $\SO_{2n+1}$, the constants $c_i\in \Qlb$ satisfy that $\sum\limits_i |c_i| \leq 2^n \sqrt{2^n n!}$ for any identification $\mathbb C \simeq \Qlb$.
\end{prop}
\begin{proof}
This is an easy corollary of \cite{L2} Theorem 5.6. Note that any irreducible unipotent representation can appear on the right hand side there for some virtual cell $\underline {c}$ by Lemma 2.22 and Lemma 5.2 of \cite{L2},  that every irreducible representation $E$ of $\mathrm W_{\SO_{2n+1}}$ has $\dim (E) \leq \sqrt{2^n n!}$, and that $d\leq n$ in the definition of the virtual cell (see Equation (2.12.1) of  \cite{L2}).

\end{proof}
\begin{rmk}
Clearly, a similar result holds for $\Sp_{2n}$, as we may argue as above or use a Lusztig map $\mathfrak L_s$ as introduced in Subsection \ref{DLc}.
\end{rmk}
~

\subsection{Unipotent irreducible representations of $\mathrm O_{2n}$}
We denote the orthogonal group of a vector space $V$ over $\mathbb F_q$ of dimension $2n$ equipped with a nondegenerate symmetric bilinear form by $\mathrm O(V)$. Note that there are two isomorphism classes of the orthogonal groups as the $2n$-dimensional space $V$ varies. Let $G_n$ be an algebraic group  isomorphic to $\mathrm O(V)$ for some $V$ as mentioned in the above.
We need a parallel result for $G_n$. We denote  the identity component of $G_n$ by $G_n^\circ$. Recall that we write $R^{G_n}_{T,\theta}$ for $Ind^{\mathrm G_n^F}_{G_n^{\circ, F }}R^{G^\circ_n}_{T,\theta}$.
\begin{prop}\label{Direp}Assume that $q \geq 2^{2n}$.
Let $\pi \in \mathcal E(G_n,1)$. Then $\pi$ is a subrepresentation of a representation $\phi$ whose character $\chi_\phi$ is of the form 
$$\chi_\phi=\sum_i c_i R_i,
$$
where each $R_i$ is of the form $ R_{T,1}^{G_n}$ for some $F$-stable maximal torus $T$ in $G_n^\circ$, the constants $c_i\in \Qlb$ satisfy that $\sum\limits_i |c_i| \leq 2^n \sqrt{2^n n!}$ for any identification $\mathbb C \simeq \Qlb$.
\end{prop}
\begin{proof}
This is a corollary of \cite{L3} Proposition 3.13. Let $\pi_0$ be a subrepresentation of the restriction of $\pi$ to $G_n^{\circ,F}$. Combining Lemma 1.19 and Lemma 3.8 of \cite{L3}, we may argue as in Proposition \ref{Birep} to get a representation $\phi_0$ of $G_n^{\circ, F}$ containing $\pi_0$ as a subrepresentation whose character $\chi_{\phi_0}$
is of the form
$$ \chi_{\phi_0}=\sum_i c_i R_{0,i},
$$
where each $R_{0,i}$  is of the form $R_{T,1}^{G_n^\circ}$ for some rational maximal torus $T$ in $G_n^\circ$, the constants $c_i\in \Qlb$ satisfy that $\sum |c_i| \leq 2^n \sqrt{2^n n!}$ for any identification $\mathbb C \simeq \Qlb$. Then we make an induction to $G_n^F$, setting $R_i=Ind^{\mathrm G_n^F}_{G_n^{\circ,F}} R_{0,i}$ and $\phi = Ind^{\mathrm G_n^F}_{G_n^{\circ,F}} \phi_0$.
\end{proof}
\subsection{Fourier-Jacobi models over finite fields}\label{fjmoff}
In this subsection, we recall the Fourier-Jacobi model of $\Sp (V) \simeq\Sp_{2n}$ over finite fields, where $V$ is a vector space of dimension $2n$ equipped with a non-degenerate symplectic form. See \cite{W} Section 1.2 or \cite{GGP1} Chapter 13 for details.

Let $W$ be a subspace of $V$, so that the restriction of the symplectic form of $V$ to $W$ is non-degenerate. Let $W^\perp$ be the orthogonal complement of $W$.
Then we have $V= W \oplus W^\perp$. 
Choose a polarization $W^\perp=X \oplus X^\vee$, so that $X$ and $X^\vee$ are maximal totally isotropic subspaces of $W^\perp$ which form linear duals of each other under the form on $V$.
Fix a basis $\{x_i\}_{1\leq i\leq l}$ of $X$. 
Let $P$ be the parabolic subgroup of $\Sp(V)$ stabilizing the flag \begin{equation}\label{flag}
0\subset \Span(\{x_1\}) \subset \ldots \subset\Span({\{x_i\}_{1\leq i\leq l}})=X \subset V
\end{equation}
The Levi subgroup of $P$ can be identified with $\mathbb G_m^{ l}\times \Sp(W)$. Let $N$ be the unipotent radical of $P$. Then we have $N=N_{\GL_l}\ltimes N_l$, here $N_{\GL_l}$ is the unipotent radical of the Borel subgroup of  $\GL(X)$ stabilizing the flag (\ref{flag}),
and $N_l$ is the unipotent radical of the parabolic subgroup  $P_l$ of $\Sp(V)$ stabilizing $X$.
 There is a homomorphism $N^F_l \to \mathcal H (W)$ invariant under the conjugation action of $N_l^F$ and $N_{\GL_l}^F$, where $\mathcal H(W)$ is the Heisenberg group of $W$. Recall that we have fixed a nontrivial character $\psi: \mathbb F_q \to \Qlb^\times$. We may view the Weil representation $\omega_\psi$ as a representation of $H:=\Sp(W)^F \ltimes N^F$ by extending $\omega_\psi$ trivially to $N_{\GL_l}$. Set $\psi_l$ to be the character of $N_{\GL_l}^F$ given by
$$
\psi_l(z)=\psi(\sum \limits^{l-1}_{i=1} z_{i,i+1}),~ z\in  N^F_{\GL_l}.
$$
Set $\nu =\omega_\psi \otimes \psi_l$, regarded as a representation of $H$.

If $l\neq 0$, the multiplicity of Fourier-Jacobi models for $\pi \in \mathcal E (\Sp(V)^F) $ and $\pi' \in \mathcal E (\Sp(W)^F)$ is
$$m_\psi^{n,l}(\pi,\pi'):= \dim \Hom_{H}(\pi,\pi' \otimes \nu).
$$
If $l=0$,  the multiplicity of Fourier-Jacobi models for $\pi \in \mathcal E (\Sp(W)^F) $ and $\pi' \in \mathcal E (\Sp(W)^F)$ is
$$m_\psi^{n,0}(\pi,\pi') := \dim \Hom_{\Sp(V)^F}(\pi ,\pi' \otimes \omega_\psi).
$$
The case for $l=0$ is also referred to as the  the basic case of  Fourier-Jacobi models.

~

\subsection{Basic case for uniform functions}
In this subsection, the reductive group  $G$ is $\Sp(V) \simeq \Sp_{2n}$. We recall Theorem 4.3 of \cite{LMS}. We need some preliminaries.

For an $F$-stable maximal torus $T$ of $G$, we define $W_G(T):= N_G(T)/T$.
Let us fix a rational Borel pair $(B_0,T_0)$ of $G$. Then $\mathrm W_G:=W_G(T_0)=N_G(T_0)/T_0$ is the Weyl group of $G$. It is well-known that there is a bijection between the $G^F$-conjugacy class of $F$-stable maximal tori of  $G$ and the set of cohomology class $H^1(F,\mathrm W_G)$ (\select{c.f.}, Corollary 1.14 of \cite{DL}). For an $F$-stable maximal torus $T$ of $G$, we denote the class of $T$ in $H^1(F,\mathrm W_G)$ by ${\rm cl}(T,G)$.

For a semisimple element $s \in G$, we denote the identity component of the centralizer $C_G(s)$ of s in $G$ by $G_s$. 

If $G_s$ is $F$-stable, we define $\mathrm W_{G_s}$ in a similar way. A map 
$$
j_{G_s}: H^1(F,\mathrm W_{G_s}) \to H^1(F,\mathrm W_G)
$$
is defined in \cite{R}, which sends the class of an $F$-stable maximal torus in $G_s$ to the class of the same torus in $G$.

For $s\in S$, let $V^s$ be the subspace of $V$ fixed by $s$.
For an $F$-stable maximal torus $S$ in $G$, let $J(S)$ be an index set for the following set of symplectic subspaces of $V$:
$$
\{V^s:s \in S\}.
$$
Then $J(S)$ is finite of cardinality $|J(S)|=2^n$. The Frobenius $F$ acts naturally on $J(S)$. For $\jmath \in J(S)$, denote by $V^\jmath $ the corresponding symplectic subspace of $V$.
Set for $\jmath \in J(S)$
$$S_\jmath = \{s\in S :V^s=V^\jmath\}.$$
For $\jmath \in J(S)^F$, we set 
$$G_\jmath := C_G(S_\jmath)=Z_\jmath \times \Sp(V^\jmath).
$$
It is clear that $Z_\jmath\subset S$ can be viewed as a maximal torus of $\Sp(V/V^\jmath)$.
For $\jmath \in J(S)^F$, let $\phi(\jmath)$ be the image of a regular element of $Z_\jmath$ under the inclusion $Z_\jmath \hookrightarrow G_\jmath$. Then $G_\jmath=G_{\phi (\jmath)}$. 
The following is Theorem 4.3 of \cite{LMS}.
 In Section 4.3.2 of \cite{LMS}, the integer $l(Z_\jmath)$ is defined for  $\jmath \in J(S)^F$. 
For each $Z_\jmath$, the quadratic character $\vartheta_{Z_\jmath}: Z^F_\jmath \to \Qlb^\times$ is defined in \cite{LMS}.

\begin{thm} \label{MF} Let the definitions be as  above. Then we have the multiplicity 

$$\langle R^G_{T, \chi} \otimes\omega_\psi^\vee, R^G_{S, \eta}\rangle_{G^F}
= \sum_{\jmath\in J(S, T)}\frac{(-1)^{\sigma (T) +\sigma(  S)+ l(Z_\jmath)}}{|W_{G_\jmath}(T_\jmath)^F| \, |W_G(S)^F|}\sum_{w\in W_G(T_\jmath)^F, \, v\in W_G(S)^F} \langle {}^w\chi_\jmath \vartheta_{Z_\jmath}, {}^v\eta\rangle_{Z_\jmath^F},$$
where $J(S, T):=\left\{\jmath\in J(S)^F : j_{G_\jmath}^{-1}({\rm cl}(T, G)) \neq \varnothing\right\}$,
the pair $(T_\jmath, \chi_\jmath)$ is any $G^F$-conjugate of $(T, \chi)$ such that 
$T_\jmath\subset G_\jmath$, and we denote $\chi_\jmath \circ \Ad(w^{-1}):Z_\jmath^F \to \Qlb^\times$ by $^w\chi_\jmath$ and  $\eta \circ \Ad(v^{-1}): Z^F_\jmath \to \Qlb^\times $ by $^v\eta$. 
\end{thm}

Since $|J(S,T)|\leq |J(S)|=2^n$ and  $|W_G(T_\jmath)^F|,| W_G(S)^F|\leq 2^n n!$, we have the following
\begin{cor} \label{bu}
$$|\langle R^G_{T, \chi} \otimes\omega_\psi^\vee, R^G_{S, \eta}\rangle_{G^F}|\leq 4^n n!.
$$
\end{cor}

\begin{rmk}
Theorem 4.3 of \cite{LMS} is valid when $G$ is a unitary group. This formula combined with Proposition \ref{Airep} will imply an upper bound for the basic case of the Fourier-Jacobi models of unitary groups over finite fields which is  independent of $q$. See Section 4.5 of \cite{LMS} for details. 

\end{rmk}
~

\subsection{Reduction to uniform characters}
Since not all $\pi \in \mathcal E(\Sp_{2n})$ are uniform, we need to reduce the estimation of $m_{\psi}^{n,0}(\pi,\pi')$ to the estimation of $m_\psi^{n,0}(R^{\Sp_{2n}}_{T,\theta},R^{\Sp_{2n}}_{S,\chi})$ in order to apply Corollary \ref{bu}. 
In this  subsection $G$ will be $\Sp_{2n}$ hence the dual group is $\SO_{2n+1}$.
 Recall the bijection $\mathfrak L_s$ defined in Section \ref{DLc}. We make it explicit in this case. 
For a semisimple element $s\in \SO_{2n+1}^F$, we know that $C_{\SO_{2n+1}}(s)$ is isomorphic to
$$G^{(1)}(s) \times G^{(2)}(s) \times G^{(3)}(s) ,$$
where 
\begin{itemize}
\item
 $G^{(1)}(s)$ is the reductive group corresponding to the eigenspaces of $s$ with eigenvalues $\neq \pm 1$. Hence the finite group $G^{(1)}(s)^F$ is isomorphic to a product of general linear groups and unitary groups (possibly defined over an extension of $\mathbb F_q$); 
\item
$G^{(2)}(s)$ is the reductive group corresponding to the eigenspace of $s$ with eigenvalue $1$. Hence it is isomorphic to $\SO_{2m_1+1}$ for some nonnegative integer $m_1$;
\item
$G^{(3)}(s)$ is  the reductive group corresponding to the eigenspace of $s$ with eigenvalue $-1$. Hence it is isomorphic to $\mathrm O(V)$ for some $2m_2$-dimensional vector space $V$ equipped with a  nondegenerate symmetric bilinear form. When $m_2=0$, the group $G^{(3)}(s)$ is the trivial group.
\end{itemize}

\begin{rmk} \label{resu}
The above statement 
claims that each factor of $G^{(1)}(s)$ is of the form $\Res_{\mathbb F_{q^{d}}/\mathbb F_q} \GL_{r} $ or $\Res_{\mathbb F_{q^{d}}/\mathbb F_q} \mathrm U_{r}  $ for some integers $r$ and $d$, where we denote the Weil restriction from $\mathbb F_{q^d}$ by  $\Res_{\mathbb F_{q^d}/\mathbb F_q}$.
To apply Corollary \ref{Airep2} and Remark \ref{Airep3} to unipotent representations of  $G^{(1)}(s)^F$, we note that the classes of unipotent representations are preserved under the natural isomorphisms $\Res_{\mathbb F_{q^{d}}/\mathbb F_q} \GL_{r}(\mathbb F_q) \simeq \GL_{r}(\mathbb F_{q^{d}})$ and $\Res_{\mathbb F_{q^{d}}/\mathbb F_q} \mathrm U_{r}(\mathbb F_q) \simeq \mathrm U_{r}(\mathbb F_{q^{d}})$. See (1.18) of \cite{L4}.
\end{rmk}

We have a bijection $\mathfrak L_s:\mathcal E(G,s) \to \mathcal E(G^{(1)}(s),1) \times 
\mathcal E(G^{(2)}(s),1) \times \mathcal E(G^{(3)}(s),1) $. Since the linear extension of  $\mathfrak L_s$ 
sends Deligne-Lusztig (virtual) characters to Deligne-Lusztig characters up to a sign and sends characters of representations  to characters of representations, we have the following.                 
\begin{prop} \label{embed}
Assume that $q\geq 2^{2n}$. Let $\pi \in \mathcal E(\Sp_{2n})$. Then $\pi$ is a subrepresentation of a representation $\gamma$ whose character $\chi_\gamma$ is of the form
$$\chi_\gamma =\sum_i c_i R_i ,
$$
where each $R_i$ is of the form $R_{T,\theta}^{\Sp_{2n}}$ for some $F$-stable maximal torus $T$ in $\Sp_{2n}$ and some 1-dim representation $\theta:T^F \to \Qlb^\times$, and the constants $c_i\in \Qlb$ satisfy that $\sum\limits_i |c_i| \leq 2^n \sqrt{2^n n!}$ for all identifications $\mathbb C \simeq \Qlb$.
\end{prop}
\begin{proof}
Suppose that $\pi \in \mathcal E(\Sp_{2n},s)$. Fix a Lusztig map $\mathfrak L_s$.
We denote $C_{G^*}(s)= G^{(1)}(s)\times G^{(2)}(s) \times G^{(3)}(s)$ by $ H$ to simplify the notation.

Step 1: We reduce the problem to the case concerning unipotent representations of $(G^{(1)}(s)\times G^{(2)}(s) \times G^{(3)}(s))^F$. If the representation $\mathfrak L_s(\pi)$ can be embedded into a representation $\gamma_1$ whose character $\chi_{\gamma_1}$ is of the form 
$$\chi_{\gamma_1} = \sum _i \hat c_{i} R_{i,1},
$$
where each $R_{i,1}$ is of the form $R_{S,1}^{ H}$ for some $F$-stable maximal torus $S$ in the identity component of $ H$, and $\sum\limits_i |\hat c_i|\leq 2^n \sqrt{2^n n!}$ for all identifications $\mathbb C \simeq \Qlb$. Then on the one hand, we have 
$$\mathfrak L_s^{-1}(\chi_{\gamma_1})= \sum_i c_i R_i , 
$$
where $R_i$ are Deligne-Lusztig characters satisfying $\mathcal L_s (R_i)= \pm R_{i,1}$ and we have $c_i =\pm \hat c_i$. On the other hand,
we write $\gamma_1= \mathfrak L_s(\pi) \oplus \mathfrak L_s(\pi)^\vee$ for some representation $\mathfrak L_s(\pi)^\vee$ of $ H^F$. We denote the representation $\mathfrak L_s^{-1}(\mathfrak L_s(\pi)^\vee)$ by $\pi^\vee$. Then $$
\mathfrak L_s^{-1}(\chi_{\gamma_1})= \chi_{\pi}+\chi_{\pi^\vee},
$$
where $\chi_\pi$ is the character of $\pi$ and $\chi_{\pi^\vee}$ is the character of $\pi^\vee$.
In particular, we see that 
$\mathfrak L_s^{-1}(\chi_{\gamma_1})$ is the character of a representation $\gamma$ of $G^F$ which contains $\pi$ as a subrepresentation. Hence $\pi$ is a subrepresentation of a representation $\gamma$ whose character $\mathfrak L_s ^{-1}(\chi_{\gamma_1})$ is of the desired form.
We have reduced the problem to find the desired $\gamma_1$ whose character $\chi_{\gamma_1} = \sum \limits_i \hat c_{i} R_{i,1}$ satisfies $\sum\limits_i |\hat c_i|\leq 2^n \sqrt {2^n n!}$.

Step 2:
Combining Corollary \ref{Airep2}, Proposition \ref{Birep} and Proposition  \ref{Direp}.
More precisely, let $\mathfrak L_s (\pi)=\pi^{(1)} \boxtimes \pi^{(2)}\boxtimes \pi^{(3)}$ be the decomposition along $ H^F=G^{(1),F}(s) \times  G^{(2),F}(s)\times G^{(3),F}(s)$. Using Remark \ref{resu}, let $\chi_{\pi^{(1)}}=\sum\limits_i c^{(1)}_i R_i$ be as in Corollary \ref{Airep2}  where $R_i$ are unipotent Deligne-Lusztig characters of $G^{(1),F}(s)$. Let $\chi_\tau=\sum\limits_j c^{(2)}_j R_j$ be as in Proposition \ref{Birep} for some representation $\tau$ containing $\pi^{(2)}$ as a subrepresentation for unipotent Deligne-Lusztig characters $R_j$ of $G^{(2),F}(s)$. Let $\chi_\iota = \sum\limits_k c^{(3)}_k R_k $ be as in Proposition \ref{Direp} for some representation $\iota$ containing $\pi^{(3)}$ as a subrepresentation for unipotent Deligne-Lusztig characters $R_k$ of $G^{(3),F}(s)$. Let $\gamma_1=\pi^{(1)}\boxtimes \tau \boxtimes \iota$ be the representation whose character is $\chi_{\pi^{(1)}}\boxtimes \chi_\tau \boxtimes \chi_\iota$. Note that for each
$(i,j,k)$, the virtual character $R_i \boxtimes R_j \boxtimes R_k$ is still a unipotent Deligne-Lusztig character. 
We denote $c_{i,j,k}:=c_i^{(1)}c_j^{(2)}c_k^{(3)}$ and $R_{i,j,k}:=R_i \boxtimes R_j \boxtimes R_k$.
To conclude,
 we can embed $\mathfrak L_s(\pi)$ into a representation $\gamma_1$ whose character $\chi_{\gamma_1}$ is of the form
 $$\chi_{\gamma_1} = \sum _{i,j,k}  c_{i,j,k} R_{i,j,k},
$$
where   $\sum\limits_{i,j,k} | c_{i,j,k}|\leq 2^n \sqrt{2^n n!}$ (as can be easily verified using the estimation in Corollary \ref{Airep2}, Proposition \ref{Birep} and Proposition \ref{Direp}). This completes the proof.
\end{proof}

\begin{rmk}
A similar and general result is proved in Lemma A.1 of \cite{S}.
\end{rmk}

\section{An Upper Bound}
In this section, we assemble all the ingredients at hand.
~

\subsection{Estimation for the basic case of Fourier-Jacobi models}
Our main theorem is the following:

\begin{thm}\label{mt}
Assume that $q\geq 2^{2n}$. Let $m_\psi^{n,0}(-,-)$ be the function defined at the end of Subsection \ref{fjmoff}.
For any $\pi,\pi'\in \mathcal E(\Sp_{2n})$, we have
$$m_\psi^{n,0}(\pi,\pi')=\langle \pi \otimes \omega_\psi^\vee ,\pi'\rangle_{\Sp_{2n}^F} \leq 2^{5n}(n!)^2.
$$
\end{thm}

\begin{proof}
The equality follows from the definition of $m_\psi^{n,0}(-,-)$.
We will abuse the notation by denoting the $\Qlb$-span of $m_\psi^{n,0}(-,-)$ again by  $m_\psi^{n,0}(-,-)$. 
Let $\gamma_\pi$ and $\gamma_{\pi'}$ be as in Proposition \ref{embed} corresponding to $\pi$ and $\pi'$ respectively.
Assume that $\chi_{\gamma_\pi}=\sum\limits_i c_i R_i$ and $\chi_{\gamma_{\pi'}}=\sum\limits_j d_j R_j$ are as in Proposition \ref{embed}.
 Then we have
$$m_\psi^{n,0}(\pi,\pi') \leq m_\psi^{n,0} (\chi_{\gamma_\pi},\chi_{\gamma_{\pi'}}) \leq (\sum_i|c_i|)\cdot (\sum_j|d_j|) \cdot 4^n n! \leq 2^{5n}(n!)^2 ,
 $$ 
where the first inequality follows from the fact that $\pi$ (resp. $\pi'$) is a subrepresentation of $\gamma_\pi$ (resp. $\gamma_{\pi'}$) and the fact that $m_\psi^{n,0}(-,-)$ takes nonnegative values for a pair of characters of representations,   the second inequality follows from Corollary \ref{bu} and the third inequality follows from the properties of $c_i$ and $d_i$ displayed in Proposition \ref{embed}.

\end{proof}
~

\subsection{Reduction to basic case}
The following can be derived from Proposition 7.3 and Proposition 7.8 of \cite{W}. 
\begin{prop}\label{rbc}
Let $s$ be a semisimple element of $\SO_{2n+1}^F=(\Sp_{2n}^*)^F$, and $s'$ be a semisimple element of $\SO_{2m+1}^F=(\Sp_{2m}^*)^F$. Let $\pi\in  \mathcal E(\Sp_{2n},s)$ and $\pi' \in \mathcal E(\Sp_{2m},s') $ with $n \geq m$. Let $P$ be an $F$-stable maximal parabolic subgroup of $\Sp_{2n}$ with Levi factor $\GL_{n-m} \times \Sp_{2m}$. Let $s_0$ be a semisimple element of $\GL_{n-m}^F$ and let $\tau \in \mathcal E(\GL_{n-m},s_0)$ be an irreducible cuspidal representation of $\GL_{n-m}^F$ which is nontrivial if $n-m=1$. Assume that $s_0$ has no common eigenvalue with $s$ and $s'$. Set $l=n-m$. Let $m_\psi^{n,l}(-,-)$ be the function defined in Subsection \ref{fjmoff}.
Then we have
$$
m_\psi^{n,l}(\pi,\pi')= \langle \pi \otimes {\omega}_\psi^\vee ,I^{\Sp_{2n}}_{P}(\tau \otimes \pi')\rangle_{\Sp_{2n}^F},
$$
where we denote the parabolic induction from $P^F$ to $\Sp_{2n}^F$ by $I^{\Sp_{2n}}_P$.
\end{prop}
\begin{rmk}\label{pu}
The linear extension of the parabolic induction (we will abuse the  notation by denoting it again by $I^{\Sp_{2n}}_P$) preserves uniform functions. More precisely, we have
$$I^{\Sp_{2n}}_P (R^L_{T,\theta})= R^{\Sp_{2n}}_{T,\theta},
$$
where $L$ is a Levi subgroup of $P$, and $T$ is a $F$-stable maximal torus in $L$ which can be viewed as an $F$-stable maximal torus in $\Sp_{2n}$ via the  inclusion $L \hookrightarrow \Sp_{2n}$.
\end{rmk}

\begin{rmk} \label{qsr} We keep the notation and the assumption as in Proposition \ref{rbc}.
Suppose that $n>m$ and $q\geq 2^{2n}$. Then for any pair $(\pi,\pi')$ where $\pi \in \mathcal E(\Sp_{2n})$ and $\pi' \in \mathcal E(\Sp_{2m})$, we can take $\tau \in \mathcal E(\GL_{n-m},s_0)$ satisfying the assumption of Proposition \ref{rbc}, so that the character of $\tau$ is of the form $\pm R^{\mathrm {GL}_{n-m}}_{T,\chi}$ for some $F$-stable maximal torus 
$T$ of $\mathrm {GL}_{n-m}$ and some character $\chi: T^F \to \Qlb^\times$.

\end{rmk}
Hence we conclude the following.

\begin{thm}\label{mt2}
Assume that $q\geq 2^{2n}$.  Let $\pi \in \mathcal E(\Sp_{2n})$ and $\pi' \in \mathcal E(\Sp_{2m})$ with $n\geq m$. Set $l=n-m$.
Let $m_\psi^{n,l}(-,-)$ be the function defined in Subsection \ref{fjmoff}.
Then we have
$$
m_\psi^{n,l}(\pi,\pi')\leq 2^{5n} (n!)^2.
$$
\end{thm}
\begin{proof} 
Let $\gamma_\pi=\sum c_i R_i^\pi$ be a representation of $\Sp_{2n}^F$ as in Proposition \ref{embed} containing $\pi$ as a subrepresentation, and  let $\gamma_{\pi'}=\sum d_j R_j^{\pi'}$ be a representation of $\Sp_{2m}^F$ as in Proposition \ref{embed} containing $\pi'$ as a subrepresentation. We have
$$m_\psi^{n,l}(\pi, \pi')\leq m_\psi^{n,l}(\gamma_\pi,\gamma_{\pi'}) \leq \sum_{i,j} |c_i d_j| |m_\psi^{n,l}(R_i^\pi,R_j^{\pi'})| \leq 2^{5n} (n!)^2.
$$
Here the first and  second inequality is proved as in  Theorem \ref{mt}, and the last inequality follows from Remark \ref{pu}, Remark \ref{qsr}, Proposition \ref{rbc} and Corollary \ref{bu} (and the fact that $m\leq n$).
\end{proof}

\begin{rmk}\label{rmmt}
Theorem \ref{mt} implies the type $B$ case of Conjecture 3.1 of \cite{HS} (\select{i.e.}, Theorem \ref{pmt}). To get a function of $n$ independent of $q$ in down to earth terms, we may use the following trivial estimation for $q \leq 2^{2n}$ and $\pi ,\pi' \in \mathcal E(\Sp_{2n})$
$$
 \langle \pi \otimes \omega^\vee_\psi, \pi' \rangle \leq \langle \pi \otimes \omega^\vee_\psi,reg \rangle 
\leq \sqrt{|\Sp_{2n}^F|} q^n \leq \sqrt{|\GL_{2n}^F|} q^{n}\leq \sqrt{q^{4n^2}}q^{n}= q^{2n^2+n}\leq 2^{4n^3+2n^2} ,
$$ where we denote the regular representation of $\Sp_{2n}^F$ by $reg$.
The desired function can be taken to be $$f(n)= \max\{2^{5n}(n!)^2,2^{4n^3+2n^2}\}.$$
 Similarly we see from Theorem \ref{mt2} that there is an upper bound $M(n,l)$ for the multiplicities  $m_\psi^{n,l} (-,-)$ which is independent of $q$.
\end{rmk}
\begin{rmk}
It is likely that Theorem \ref{mt} holds without restrictions on $q$, see Section 3.19 of \cite{L3}.
\end{rmk}

\end{document}